\documentclass{article}

\usepackage{amsfonts}
\usepackage{amsmath}
\usepackage{amssymb}
\usepackage{amsthm}
\usepackage{spconf}
\usepackage{etoolbox}
\usepackage{graphicx}
\usepackage{mathtools}
\usepackage{subcaption}

\makeatletter
\g@addto@macro\th@plain{\thm@headpunct{\ }}
\makeatother
\makeatletter
\g@addto@macro\th@definition{\thm@headpunct{\ }}
\makeatother
\makeatletter
\def\maketag@@@#1{\hbox{\m@th\normalfont\normalsize#1}}
\makeatother

\def\mllap{\mathpalette\mllapinternal}
\def\mllapinternal#1#2{\llap{$\mathsurround=0pt#1{#2}$}}

\makeatletter
\newcommand{\pushright}[1]{\ifmeasuring@#1\else\omit\hfill$\displaystyle#1$\fi\ignorespaces}
\newcommand{\pushleft}[1]{\ifmeasuring@#1\else\omit$\displaystyle#1$\hfill\fi\ignorespaces}
\makeatother

\csundef{proof}
\csundef{endproof}
\theoremstyle{plain}\newtheorem{theorem}{Theorem}
\theoremstyle{plain}\newtheorem{corollary}{Corollary}
\theoremstyle{definition}
\theoremstyle{definition}

\title{Spectral Statistics of Lattice Graph\\ Structured, Non-uniform Percolations}
\name{Stephen Kruzick$^1$ and Jos\'{e} M. F. Moura$^2$ \thanks{This work was supported under NSF grant {\#}CCF1513936.}\thanks{Email: $^1$skruzick@andrew.cmu.edu, $^2$moura@ece.cmu.edu}}
\address{Carnegie Mellon University, Department of Electrical Engineering\\5000  Forbes Avenue, Pittsburgh, PA 15213}

\begin{document}
\maketitle

\begin{abstract}

%Design of filters for digital signal processing over networks benefits from knowledge of the spectral decomposition of matrices that encode graphs, such as the adjacency matrix and the Laplacian matrix, used to define the shift operator.  For shift matrices with real eigenvalues, which arise for symmetric graphs, the empirical spectral distribution captures the eigenvalue locations.  Under realistic circumstances, stochastic influences often affect the network structure and, consequently, the shift matrix empirical spectral distribution.  Nevertheless, deterministic functions may often be found to approximate the asymptotic behavior of empirical spectral distributions of random matrices.  This paper uses stochastic canonical equation methods developed by Girko to derive such deterministic equivalent distributions for the empirical spectral distributions of random graphs formed by structured, non-uniform percolation of a $D$-dimensional lattice supergraph.  Included simulations demonstrate the results for sample parameters.

Design of filters for graph signal processing benefits from knowledge of the spectral decomposition of matrices that encode graphs, such as the adjacency matrix and the Laplacian matrix, used to define the shift operator.  For shift matrices with real eigenvalues, which arise for symmetric graphs, the empirical spectral distribution captures the eigenvalue locations.  Under realistic circumstances, stochastic influences often affect the network structure and, consequently, the shift matrix empirical spectral distribution.  Nevertheless, deterministic functions may often be found to approximate the asymptotic behavior of empirical spectral distributions of random matrices.  This paper uses stochastic canonical equation methods developed by Girko to derive such deterministic equivalent distributions for the empirical spectral distributions of random graphs formed by structured, non-uniform percolation of a $D$-dimensional lattice supergraph.  Included simulations demonstrate the results for sample parameters.

\end{abstract}


%The plan
%	Non-uniform percolation
%	New simulations
%	blah blah blah

\section{Introduction}

Modern technological advances have produced a world of people, devices, and systems that are increasingly connected, often in intricate ways that are best described by complex networks.  In network science graphs capture for example relations among individuals in a social network context. In data science, graphs represent dependencies among streams of data generated by different sources or agents.  Such networks are frequently large, and it may be desirable to model them as random variables due to uncertainty or inherent stochastic influences in their structures.  When studying the properties of the matrices that encode the graph structure of these networks, such as the graph adjacency matrix and the graph Laplacian, spectral decompositions are often invoked.  Linear shift-invariant filtering as defined in graph signal processing represents an example application in which such eigenvalue information would be useful.  In signal processing on graphs, a matrix $W$ related to the graph structure, such as the graph adjacency matrix or graph Laplacian, defines the shift operator \cite{ASan1}\cite{DShu1}.  Filters manifest as polynomial functions $P\left(W\right)$ in the shift operator \cite{ASan1}, and decomposition of a signal defined on the nodes according to a basis of eigenvectors of $W$ play the role of the Fourier Transform \cite{ASan2}.  Because the eigenvalues of the row-normalized adjacency matrix and row-normalized Laplacian matrix are closely related to a measure of signal complexity known as total variation, eigenvalues of the shift matrix can be interpreted as frequencies \cite{DShu1}\cite{ASan2}\cite{ASan3}\cite{SChe1}.  If $W$ is a diagonalizable matrix, then $P\left(W\right)$ is simultaneously diagonalizable with $W$, so the frequency response to an eigenvector $\mathbf{v}$ where $W\mathbf{v}=\lambda\mathbf{v}$ is $P(\lambda)$ \cite{ASan3}\cite{ASan4}.  Hence, information concerning the shift matrix eigenvalues is critical to filter design.  For instance, knowledge of the eigenvalues can lead to polynomial filters that accelerate the distributed average consensus algorithm \cite{EKok1} or, for large filter degrees and completely known eigenvalues, can even lead to polynomial filters that achieve consensus in finite time \cite{ASan4}.  This paper examines the adjacency matrix eigenvalues for a particular random network model as the size of the network grows through the asymptotic behavior of the empirical spectral distribution of the adjacency matrix, a function that counts the fraction of eigenvalues of a Hermitian matrix on the interval $\left(-\infty,x\right]$ \cite{VGir1}\cite{RCou1}.

Specifically, this paper examines random graphs formed by including each link of a $D$-dimensional lattice according to independent Bernouli trails with inclusion probability depending on the dimension of the lattice to which the link belongs, a non-uniform Bernouli link-percolation model \cite{GGri1}.  In this context, a $D$-dimensional lattice graph has nodes associated with $D$-tuples and has links between nodes if those nodes correspond to tuples that differ by exactly one symbol, generalizing the definition of a cubic lattice found in \cite{RLas1} for example.  These lattice graphs, which have at most $2^D$ adjacency matrix eigenvalues, represent good candidates to examine because the number of expected adjacency matrix eigenvalues for the resulting percolation model depends only on the number of dimensions $D$ and not on the lattice size \cite{RLas1}.  The asymptotic behavior of empirical spectral distributions is sometimes characterizable, as in the well known case of Wigner matrices with the semicircular law \cite{EWig1} and, closely related, in Erd\"{o}s-R\'{e}nyi model adjacency matrices \cite{XDin1}.  For the matrices in this paper, sequences of deterministic functions that asymptotically approximate the empirical spectral distribution are computed using the stochastic canonical equations tools developed by Girko, which allow analysis when symmetric matrix entries are independent, except when determined by symmetry, but not necessarily identically distributed \cite{VGir1}.  These tools were used by others to analyze the empirical spectral distribution of a different type of random network model known as stochastic block models in \cite{KAvr1}, that leads to a different system of equations and solution form when the analysis tools are applied.

%Linear shift-invariant filtering as defined in digital signal processing on graphs represents an example application in which such eigenvalue information would be useful.  In signal processing on graphs, a matrix $W$ related to the graph structure, such as the graph adjacency matrix or graph Laplacian, defines the shift operator \cite{ASan1}\cite{DShu1}.  Filters manifest as polynomial functions $P\left(W\right)$ in the shift operator \cite{ASan1}, and decomposition of a signal defined on the nodes according to a basis of eigenvectors of $W$ imitates the Fourier Transform \cite{ASan2}.  Because the eigenvalues of the row-normalized adjacency matrix and row-normalized Laplacian matrix are closely related to a measure of signal complexity known as total variation, eigenvalues of the shift matrix can be interpreted as frequencies \cite{DShu1}\cite{ASan2}\cite{ASan3}\cite{SChe1}.  If $W$ is a diagonalizable matrix, then $P\left(W\right)$ is simultaneously diagonalizable with $W$, so the frequency response to an eigenvector $v$ where $W\mathbf{v}=\lambda\mathbf{v}$ is $P(\lambda)$ \cite{ASan3}\cite{ASan4}.  Hence, information concerning the shift matrix eigenvalues is critical to filter design.  For instance, knowledge of the eigenvalues can lead to polynomial filters that accelerate the distributed average consensus algorithm \cite{EKok1} or, for large filter degrees and completely known eigenvalues, can even lead to polynomial filters that achieve consensus in finite time \cite{ASan4}.

Section \ref{Background} introduces spectral statistics concepts and an important theorem, used for the computation of the main results.  This theorem provides a method to compute deterministic functions approximating the empirical spectral distributions.  Subsequently, Section \ref{MainResults} derives results describing the deterministic equivalents for the empirical spectral distribution of non-uniform percolations of lattice graphs with arbitrary parameters.  Proof of these results is omitted for brevity but may be found by referring to \cite{SKru1}, a paper by the authors which contains these results formulated for uniform percolations of lattice graphs.  Finally, Section \ref{Conclusion} presents a few concluding remarks.

\vspace{-8pt}
\section{Background}\label{Background}

Given a $N\times N$ Hermitian matrix $W_N$ with eigenvalues ordered such that $\lambda_i\left(W_N\right)\leq\lambda_j\left(W_N\right)$ for $1\leq i < j \leq N$, the empirical spectral distribution of $W_N$ \cite{RCou1}
\begin{equation}
F_{W_N}(x)=\frac{1}{N}\sum_{i=1}^N \chi\left(\lambda_i\left(W_N\right)\leq x\right)
\end{equation}
counts the number of eigenvalues on $\left(-\infty,x\right]$.  The corresponding empirical spectral density function of $W_N$ \cite{RCou1}
\begin{equation}
f_{W}\left(x\right)=\frac{1}{N}\sum_{i=1}^N \delta\left(x-\lambda_i\left(W_N\right)\right)
\end{equation}
indicates the locations of the eigenvalues.  When $W_N$ is a random matrix, $F_{W_N}$ and $f_{W_N}$ are function-valued random variables.  The methods used to analyze the empirical spectral distribution often rely on the Stieltjes transform
\begin{equation}
\begin{aligned}
S_F\left(z\right)&=\int_{-\infty}^{\phantom{-}\infty}\frac{1}{x-z}dF\left(x\right), \quad \operatorname{Im}\left\{z\right\}\neq 0\\
&=\frac{1}{N}\operatorname{tr}\left(\left(W_N-zI_N\right)^{-1}\right), \enskip \operatorname{Im}\left\{z\right\}\neq 0
\end{aligned}
\end{equation}
  The values of $F_{W_N}$ and $f_{W_N}$ can be found by inverting the Stieltjes transform
\begin{gather}
F\left(x\right)=\lim_{\epsilon\rightarrow 0^+}\frac{1}{\pi}\int_{-\infty}^{x}\operatorname{Im}\left\{S_F\left(\lambda+\epsilon i\right)\right\}d\lambda\\ \displaybreak
f\left(x\right)=\lim_{\epsilon\rightarrow 0^+}\frac{1}{\pi}\operatorname{Im}\left\{S_F\left(x+\epsilon i\right)\right\}
\end{gather}
For a sequence of random Hermitian matrices $W_N$ indexed by $N$ and a sequence of functionals $g_N$, a sequence of deterministic matrices such that 
$
\lim_{N\rightarrow\infty}{\left(g_N\left(W_N\right)-g_N\left(W_{N}^{\circ}\right)\right)}=0
$
is known as a deterministic equivalent, with $g_N\left(W_N^\circ\right)$ also called a deterministic equivalent of $g_N\left(W_N\right)$ \cite{RCou1}.
Theorem~\ref{GirkoK01Thm}, which is the primary theorem used in the computations of this paper, provides a method of computing a deterministic equivalent for the empirical spectral distribution of random symmetric matrices with independent entries in the upper triangular region and other entries determined by symmetry \cite{VGir1}.

\begin{theorem}[Girko's K1 Equation \cite{VGir1}]\label{GirkoK01Thm}

Consider a family of symmetric matrix valued random variables $W_N$ indexed by size $N$ such that $W_N$ is an $N\times N$ symmetric matrix in which the entries on the upper triangular region are independent.  That is, $\left\{\left(W_{N}\right)_{ij}|1\leq i\leq j\leq N\right\}$ are independent with $\left(W_{N}\right)_{ji}=\left(W_{N}\right)_{ij}$.  Let $W_N$ have expectation $B_N=\operatorname{E}\left[W_N\right]$ and centralization $H_N=W_N-E\left[W_N\right]$ such that the following three conditions hold.  Note that in order to avoid cumbersome indexing, the index $N$ will henceforth be omitted from most expressions involving $W_N$, $B_N$, and $H_N$.
\begin{gather}
\sup_{N}{\max_{i}{\sum_{j=1}^{N}{\left|B_{ij}\right|}}}<\infty \label{GirkoCond1} \\
\sup_{N}{\max_{i}{\sum_{j=1}^{N}{\operatorname{E}\left[H_{ij}^2\right]}}}<\infty \label{GirkoCond2}\\
\begin{aligned}
\lim_{N\rightarrow\infty}{\max_{i}{\sum_{j=1}^{N}{\operatorname{E}\left[H_{ij}^2\chi\left(\left|H_{ij}\right|>\tau\right)\right]}}}
=0~~\forall~~\tau>0
\end{aligned}\label{GirkoCond3}
\end{gather}
Then for almost all $x$,
\begin{equation}
\lim_{N\rightarrow\infty} \left|F_{W_N}\left(x\right)-F_N\left(x\right)\right|=0
\end{equation}
almost surely, where $F_N$ is the distribution with Stieltjes transform
\begin{equation}
\begin{aligned}\label{STEq}
S_{F_N}(z)=\frac{1}{N}\sum_{k=1}^{N}C_{kk}(z),\quad \operatorname{Im}\left\{z\right\}\neq 0
\end{aligned}
\end{equation}
and the analytic functions $C_{kk}\left(z\right)$ satisfy the canonical system of equations
\begingroup
\thinmuskip=\muexpr\thinmuskip*1/32\relax
\medmuskip=\muexpr\medmuskip*1/32\relax
\begin{equation}\label{GirkoK01Eq}
	\ \mathclap{\hspace{-13pt}C_{kk}\left(z\right)\hspace{-2pt}=\hspace{-2pt}\left[\hspace{-2pt}\left(\hspace{-4pt}B-zI-\left(\hspace{-4pt}\delta_{lj}\sum_{s=1}^{N}{C_{ss}(z)\operatorname{E}\left[H_{js}^2\right]}\hspace{-3pt}\right)_{\hspace{-3.5pt}l,j=1}^{\hspace{-3pt}l,j=N}\hspace{-2pt}\right)^{\hspace{-6pt}-1}\hspace{-1pt}\right]_{\hspace{-3pt}kk}}
\end{equation}
\endgroup
for $k=1,\ldots,N$.  Note that the notation $\left(\cdot\right)_{l,j=1}^{l,j=N}$ indicates a matrix built from the parameterized contents of the parentheses, such that $X=\left(X_{ij}\right)_{l,j=1}^{l,j=N}$, and $\delta_{lj}$ is the Kronecker delta function.  There exists a unique solution $C_{kk}(z)$ for $k=1,\ldots,N$ to the canonical system of equations \eqref{GirkoK01Eq} among %the class
$L=\left\{X(z)\in \mathbb{C} \mid X\left(z\right)\textrm{ analytic},~\operatorname{Im}\left\{z\right\}\operatorname{Im}\left\{X\left(z\right)\right\}>0\right\}$.  Furthermore, if
\begin{equation}
\inf_{i,j} {N\operatorname{E}\left[H_{ij}^2\right]}\geq c >0, \label{GirkoCond4}
\end{equation}
then
\begin{equation}
\lim_{N\rightarrow\infty} {\sup_{x}{\left|F_{W_N}\left(x\right)-F_N\left(x\right)\right|}}=0
\end{equation}
almost surely, where $F_N$ is defined as above.

\end{theorem}

\section{Main Results}\label{MainResults}

%This results of the paper contained in this section pertain to computing deterministic approximations to the empirical spectral distribution of Bernoulli link-percolation models with lattice supergraphs.  Definitions of lattice graphs vary in literature, so the relevant definition will be made precise here.  In a $D$-dimensional lattice graph with size $M_d$ along the $d$th dimension, the $\left|\mathcal{V}\right|=N=\prod_{d=1}^{D}{M_d}$ nodes are identified with the ordered $D$-tuples, where the $d$th entry has $M_d$ possible symbols.  A link connects two nodes if the corresponding $D$-tuples differ by exactly one symbol \cite{RLas1}.  Note that any integer $1\leq x \leq N$ can be written in a mixed-radix system as
This section computes deterministic approximations to the empirical spectral distribution of Bernoulli link-percolation models with lattice supergraphs.  Definitions of lattice graphs vary in the literature, so the relevant definition will be made precise here.  In a $D$-dimensional lattice graph with size $M_d$ along the $d$th dimension, the $\left|\mathcal{V}\right|=N=\prod_{d=1}^{D}{M_d}$ nodes are identified with the ordered $D$-tuples, where the $d$th entry has $M_d$ possible symbols.  A link connects two nodes if the corresponding $D$-tuples differ by exactly one symbol \cite{RLas1}.  Note that any integer $1\leq x \leq N$ can be written in a mixed-radix system as
\begin{equation}
x=1+\sum_{d=1}^{D}{\beta\left(x,d\right)\left(\prod_{j=1}^{d-1}{M_j}\right)}
\end{equation}
for $0\leq \beta\left(x,d\right)\leq M_d-1$.  Collecting the digits into a vector $\beta\left(x\right)$, the adjacency matrix of the lattice graph adjacency matrix may be written as
\begin{equation}
A_{ij}\left(\mathcal{G}_{\mathrm{lat}}\right)=\left\{\begin{array}{cc} 1 & \left\|\beta\left(i\right)-\beta\left(j\right)\right\|_0=1 \\ 0 & \textrm{otherwise} \end{array}\right.
\end{equation}
or, in terms of Kronecker products, may be written as
\begin{equation}
A\left(\mathcal{G}_{\mathrm{lat}}\right)=\sum_{j=1}^{D}{\bigotimes_{d=1}^{D}{X_{dj}}}, \enskip X_{dj}=\left\{\begin{array}{ll} K_{M_d} & j=d \\ {\phantom{K}\mllap{I}}_{M_d} & j\neq d  \end{array}\right.
\end{equation}
where $K_{M_d}$ is the complete graph on $M_d$ nodes.  The random graph model $\mathcal{G}_{\mathrm{perc}}\left(\mathcal{G}_{\mathrm{lat}},\left\{p\right\}_{d=1}^{d=D}\right)
$ under consideration starts with a $D$-dimensional lattice supergraph of given size parameters $\mathcal{G}_{\mathrm{lat}}\left(\vphantom{M_d}\right.\hspace{-3pt}\left\{M_d\right\}_{d=1}^{d=D}\hspace{-3pt}\left.\vphantom{M_d}\right)$ and includes each link of the supergraph according to an independent Bernoulli trial with probability $p_d$ depending on the lattice dimension along which the supergraph link exists, forming a non-uniform Bernoulli percolation model \cite{GGri1}.  One can verify that the scaled adjacency matrix $W\left(\mathcal{G}_{\mathrm{perc}}\right)=\frac{1}{\gamma}A\left(\mathcal{G}_{\mathrm{perc}}\right)$, where $\gamma$ is the expected node degree, satisfies the conditions for application of Theorem \ref{GirkoK01Thm}.

Theorem \ref{Thm2} derives the form \eqref{SolutionForm} of the solution to equation \eqref{GirkoK01Eq} for the scaled adjacency matrix $W\left(\mathcal{G}_{\mathrm{perc}}\right)$, and Corollary \ref{Cor1} obtains a system of equations \eqref{CoeffSystem} that describe the parameters of the solution \eqref{SolutionForm}.  This can then be used to find a deterministic equivalent for the empirical spectral distribution of $W\left(\mathcal{G}_{\mathrm{perc}}\right)$.  Note that supporting proofs of these statements have been omitted for brevity.  Corresponding proofs for the uniform percolation case $p_d=p$ for $d=1,\ldots,D$, which apply to these statements with subtle modification, can be found in \cite{SKru1}.  The proof of Theorem \ref{Thm2} relies on symmetry and Theorem \ref{GirkoK01Thm}.  Corollary \ref{Cor1} follows by simultaneous diagonalizability of terms in a matrix equation.

%\newpage
\begin{theorem}[Solution Form for $D$-Lattice Percolation]\label{Thm2}\ \\
Consider the $D$-dimensional lattice graph $\mathcal{G}_{\mathrm{lat}}$ with $N=\prod_{d=1}^{D}{M_d}$ nodes in which the $d$th dimension of the lattice has size $M_d$ for $d=1,\ldots, D$ such that the adjacency matrix is
\begingroup
\thinmuskip=\muexpr\thinmuskip*1/8\relax
\medmuskip=\muexpr\medmuskip*1/8\relax
\begin{gather}
A\left(\mathcal{G}_{\mathrm{lat}}\right)=\sum_{j=1}^{D}{\bigotimes_{d=1}^{D}{X_{dj}}}, \enskip X_{dj}=\left\{\begin{array}{ll} K_{M_d} & j=d \\ {\phantom{K}\mllap{I}}_{M_d} & j\neq d  \end{array}\right..
\end{gather}
\endgroup
Form a random graph $\mathcal{G}_{\mathrm{perc}}\left(\mathcal{G}_{\mathrm{lat}},\left\{p_d\right\}_{d=1}^{d=D}\right)$ by independently including each link of $\mathcal{G}_{\mathrm{lat}}$ along lattice dimension $d$ with probability $p_d$.  Denote the corresponding random scaled adjacency matrix $W=\frac{1}{\gamma}A\left(\mathcal{G}_{\mathrm{perc}}\right)$, expectation $B=\operatorname{E}\left[W\left(\mathcal{G}_{\mathrm{perc}}\right)\right]$, and centralization $H\left(\mathcal{G}_{\mathrm{perc}}\right)=W\left(\mathcal{G}_{\mathrm{perc}}\right)-\operatorname{E}\left[W\left(\mathcal{G}_{\mathrm{perc}}\right)\right]$ where $\gamma$ is the expected node degree
\begin{equation}
\gamma=\gamma\left(\left\{M_d\right\}_{d=1}^{d=D},\left\{p_d\right\}_{d=1}^{d=D}\right)=\sum_{d=1}^{D}p_d\left(M_d-1\right).
\end{equation}
Let $C_{kk}\left(z\right)$ for $k=1,\ldots,N$ be the unique solution to the system of equations \eqref{GirkoK01Eq} among the class $L$ guaranteed to exist by Theorem \ref{GirkoK01Thm}, and write
\begingroup
\thinmuskip=\muexpr\thinmuskip*1/32\relax
\medmuskip=\muexpr\medmuskip*1/32\relax
\begin{equation}\label{GirkoK01Sys2}
\begin{aligned}
C\left(z\right)\hspace{-2pt}=\hspace{-2pt}\left(B-zI_N-\left(\delta_{lj}\sum_{s=1}^{N}{C_{ss}\left(z\right)\operatorname{E}\left[H_{js}^{2}\right]}\right)_{l,j=1}^{l,j=N}\right)^{\hspace{-4pt}-1}\hspace{-4pt}.
\end{aligned}
\end{equation}
\endgroup
Note that $C_{kk}\left(z\right)$ is the $k$th diagonal entry of $C\left(z\right)$ and that uniqueness of $C_{kk}\left(z\right)$ implies uniqueness of $C\left(z\right)$.  For some values of $\alpha_{i_1,\ldots,i_D}\left(z\right)$ for $i_1,\ldots,i_D=0,1$
\begin{equation}
\begin{gathered}
C\left(z\right)=\sum_{\mathclap{i_1,\ldots,i_D=0}}^{1}{ \alpha_{i_1,\ldots,i_D}\left(z\right)\bigotimes_{d=1}^{D}{Y_{di_d}}}, \\ Y_{di_d}=\left\{\begin{array}{ll} {\phantom{K}\mllap{K}}_{M_d} & i_d=0 \\{\phantom{K}\mllap{I}}_{M_d} & i_d=1  \end{array}\right..
\end{gathered}\label{SolutionForm}
\end{equation}

\end{theorem}

\begin{corollary}[Stieltjes Transform for $D$-Lattice Percolation]\ \\
The Stieltjes transform of the deterministic equivalent distribution function $F_n$ specified in Theorem \ref{GirkoK01Thm} for the empirical spectral distribution of $W\left(\mathcal{G}_{\mathrm{perc}}\right)$ is given by
\begin{equation}
S_{F_N}(z)=\alpha_{1,\ldots,1}\left(z\right), \quad \operatorname{Im}\left(z\right)\neq 0
\end{equation}
where the $2^D$ complex valued variables $\alpha_{i_1,\ldots,i_D}(z)$ for $i_1,\ldots,i_D=0,1$ solve the system
\begingroup
%\small
\thinmuskip=\muexpr\thinmuskip*1/32\relax
\medmuskip=\muexpr\medmuskip*1/32\relax
\begin{equation}
\begin{aligned}
%\pushleft{\hspace{1em}\sum_{\mathclap{i_1,\ldots,i_D=0}}^1{\alpha_{i_1,\ldots,i_D}(z)\prod_{d=1}^{D}{\lambda_{di_d}\left(j_d\right)}}=}\\
%\pushleft{\frac{1}{\frac{p}{\gamma}\left(\sum\limits_{d=1}^{D}\lambda_{d0}(j_d)\right)-z-\frac{p\left(1-p\right)}{\gamma^2}\left(\sum\limits_{d=1}^{D}{\left(M_d-1\right)}\right)\alpha_{1,\ldots,1}\left(z\right)}}
\hspace{1em}\sum_{\mathclap{i_1,\ldots,i_D=0}}^1\alpha_{i_1,\ldots,i_D}&(z)\prod_{d=1}^{D}{\lambda_{di_d}\left(j_d\right)}=\\
&\hspace{-3em}\left(\frac{1}{\gamma}\left(\sum\limits_{d=1}^{D}p_d\lambda_{d0}(j_d)\right)-z-\ldots\right.
\\&\hspace{-3em}\left.\ldots\frac{1}{\gamma^2}\left(\sum\limits_{d=1}^{D}{p_d\left(1-p_d\right)\left(M_d-1\right)}\right)\alpha_{1,\ldots,1}\left(z\right)\right)^{-1}
\end{aligned}
\label{CoeffSystem}
\end{equation}
\endgroup
of $2^D$ rational equations for $j_1,\ldots,j_D=0,1$ where
\begin{equation}
\lambda_{di_d}(j_d)=\left\{\begin{array}{ll} M_d-1 &  i_d=0,j_d=0 \\ \phantom{M_d-1}\mllap{-1} & i_d=0,j_d=1 \\ \phantom{M_d-1}\mllap{1} & i_d=1\phantom{j_d=0}\end{array}\right..
\end{equation}
\label{Cor1}
\end{corollary}

\begin{figure}[t]
\centering
\begin{subfigure}{.45\textwidth}
\begin{tabular}{cccc}
\multicolumn{2}{c}{\includegraphics[width=.45\textwidth]{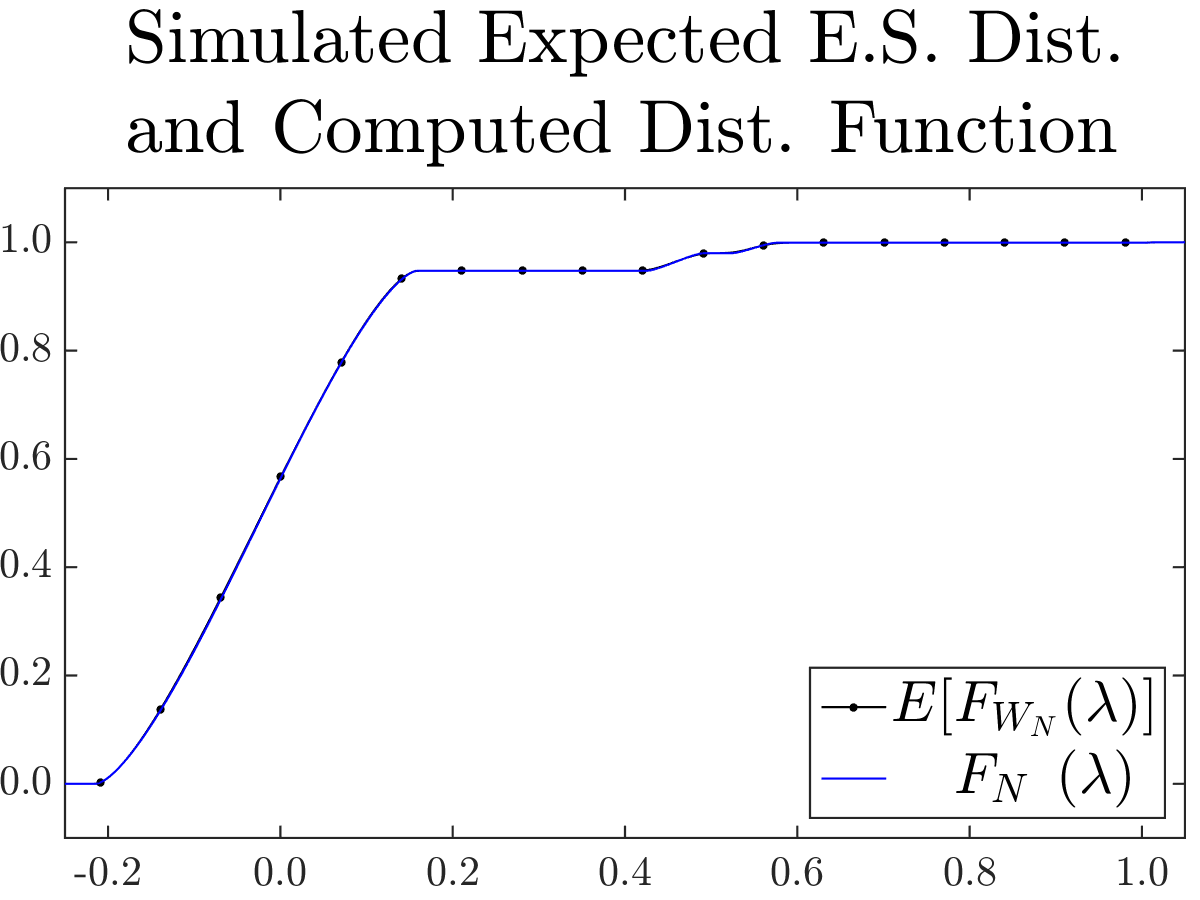}} & \multicolumn{2}{c}{\includegraphics[width=.45\textwidth]{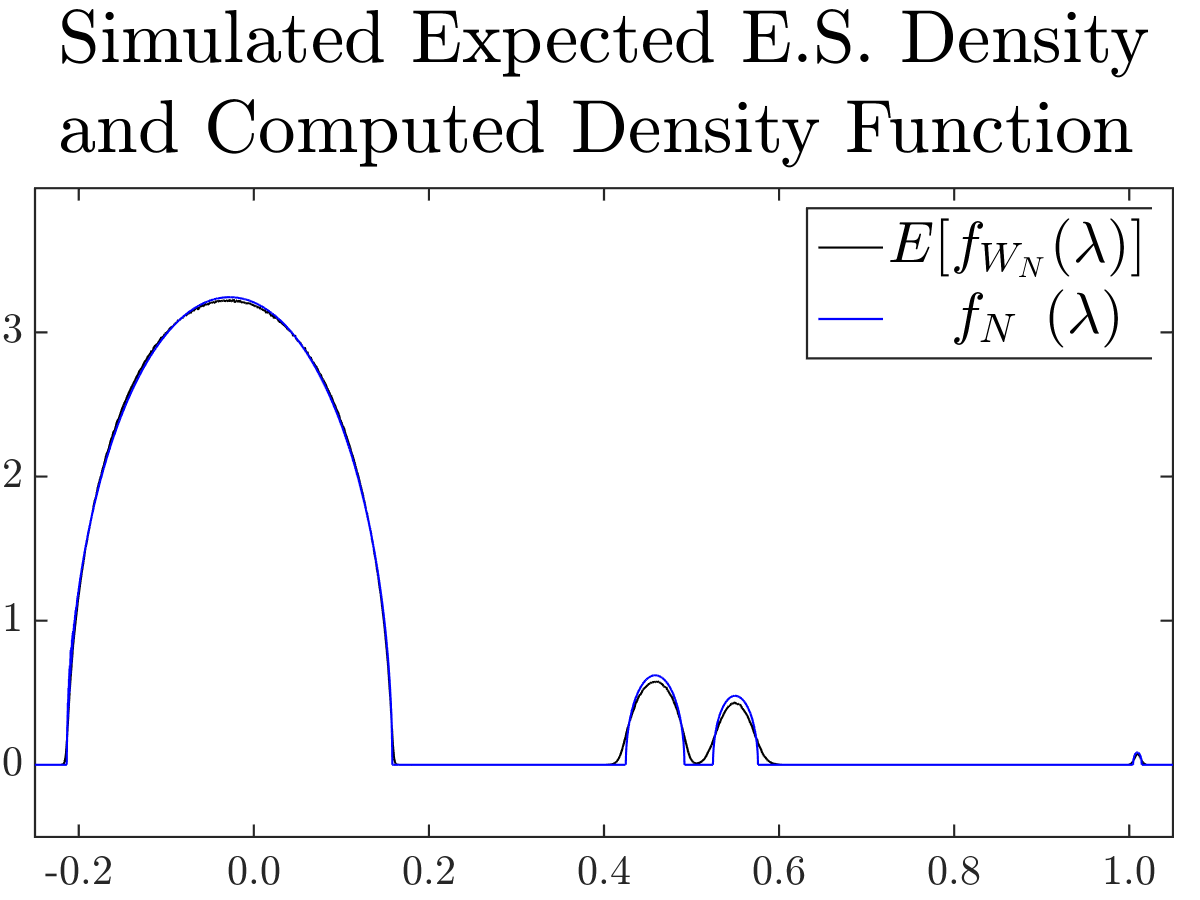}}%\\
%\includegraphics[width=.2\textwidth]{Figures/Example4_ESDens_PartA.png} & \includegraphics[width=.2\textwidth]{Figures/Example4_ESDens_PartB.png} & \includegraphics[width=.2\textwidth]{Figures/Example4_ESDens_PartC.png} & \includegraphics[width=.2\textwidth]{Figures/Example4_ESDens_PartD.png}
%\\ \multicolumn{2}{c}{\includegraphics[width=.45\textwidth]{Figures/Example4_ESDist.png}} & \multicolumn{2}{c}{\includegraphics[width=.45\textwidth]{Figures/Example4_ESDens.png}}
\end{tabular}
\caption{Comparison for 2-dimensional lattice supergraph with dimensions $(30,50)$ and percolation probabilities $(.7,.5)$.}
\label{FigEx1}
\end{subfigure}
\par\bigskip
\begin{subfigure}{.45\textwidth}
\begin{tabular}{cccc}
\multicolumn{2}{c}{\includegraphics[width=.45\textwidth]{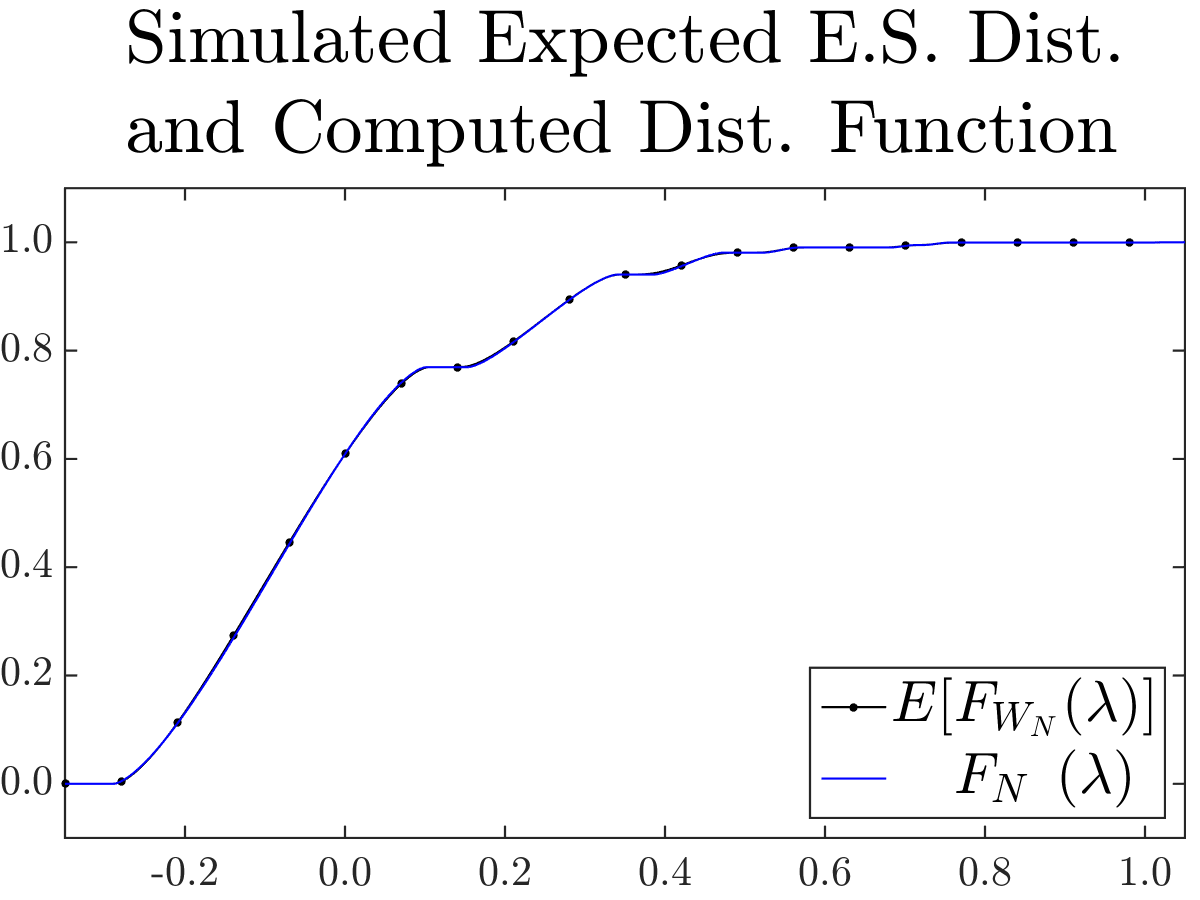}} & \multicolumn{2}{c}{\includegraphics[width=.45\textwidth]{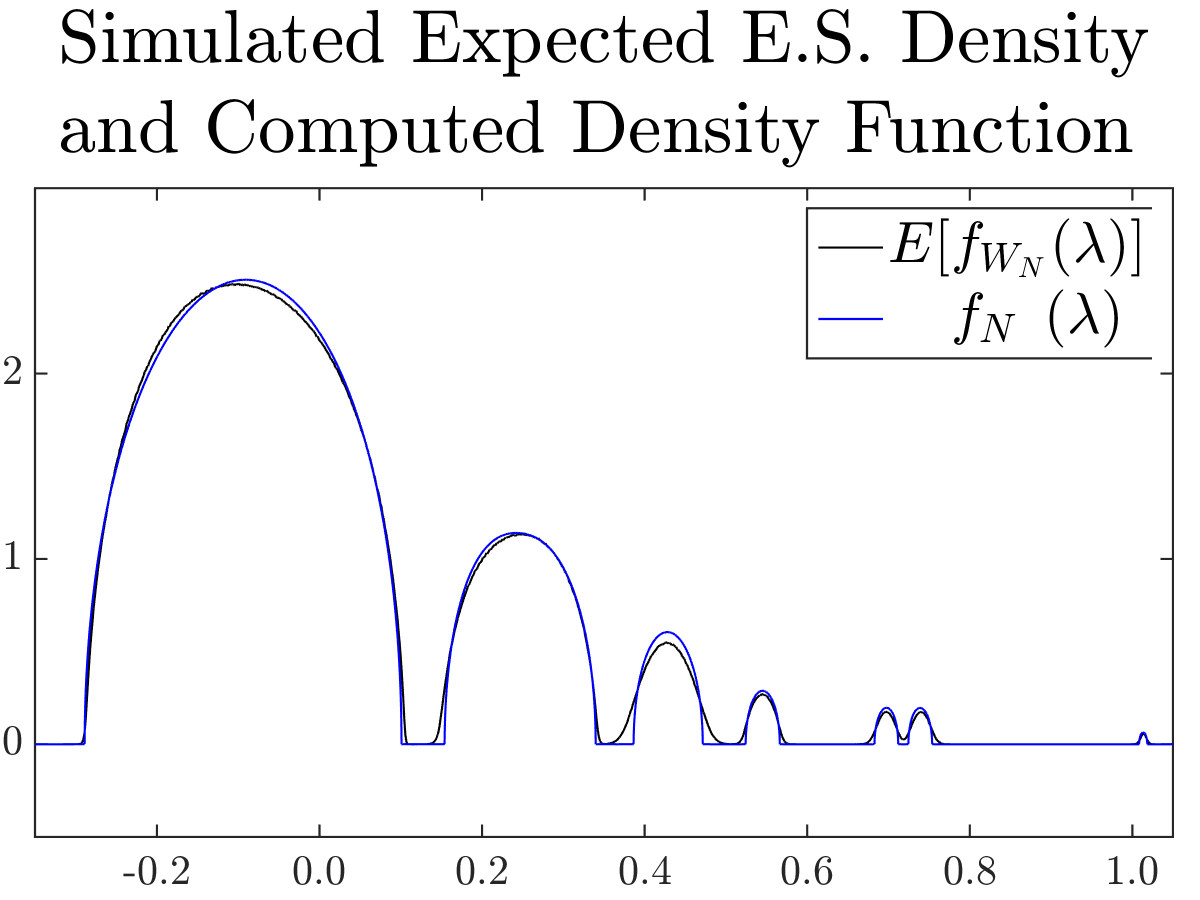}}%\\
%\includegraphics[width=.2\textwidth]{Figures/Example4_ESDens_PartA.png} & \includegraphics[width=.2\textwidth]{Figures/Example4_ESDens_PartB.png} & \includegraphics[width=.2\textwidth]{Figures/Example4_ESDens_PartC.png} & \includegraphics[width=.2\textwidth]{Figures/Example4_ESDens_PartD.png}
%\\ \multicolumn{2}{c}{\includegraphics[width=.45\textwidth]{Figures/Example5_ESDist.png}} & \multicolumn{2}{c}{\includegraphics[width=.45\textwidth]{Figures/Example5_ESDens.png}}
\end{tabular}
\caption{Comparison for 3-dimensional lattice supergraph with dimensions $(10,10,20)$ and percolation probabilities $(.8,.7,.6)$.}
\label{FigEx2}
\end{subfigure}
\caption{The above plots compare the expected empirical spectral distribution (left, black) of the scaled adjacency matrix to the computed deterministic distribution function (left, blue), with corresponding density functions (right, black and blue, respectively) also displayed, for different sample parameters.}
\label{FigEx}
\end{figure}

Solving the system of equations \eqref{CoeffSystem} for $\alpha_{1,\ldots,1}\left(z\right)$ results in an equation for which $\alpha_{1,\ldots,1}\left(z\right)$ is the unique solution with $\operatorname{Im}\left\{z\right\}\operatorname{Im}\left\{\alpha_{1,\ldots,1}\left(z\right)\right\}>0$ for $\operatorname{Im}\left\{z\right\}\neq 0$.  Consequently, it can be found for a given $z$ via zero finding methods.  Subsequently, the empirical spectral distribution can be computed by inverting the Stieltjes transform \cite{RCou1}.  Figure \ref{FigEx} shows a comparison of the computed deterministic distributions and densities to simulated expected empirical spectral distributions and densities for selected lattice parameters and percolation parameters listed in the caption.  Note that as each lattice dimension grows without bound in size, the area of the largest region of the density function asymptotically approaches totality and all other regions diminish.  Consequently, Theorem \ref{GirkoK01Thm} does not give any guarantees about the smaller regions of the density function, but a good approximation seems to be achieved.  Finally, Theorem \ref{Thm3} shows that, asymptotically, the row-normalized adjacency matrix has a similar empirical spectral distribution to that of the scaled adjacency matrix.  Hence, the computed deterministic distributions contain useful information about the scaled adjacency matrix as well.  As previously, the supporting proof can be found in \cite{SKru1}.

%\newpage
\begin{theorem}[Normalized Adjacency Matrix E.S.D.]\ \\
Let $W\left(\mathcal{G}_{\mathrm{perc}}\right)=\frac{1}{\gamma}A\left(\mathcal{G}_{\mathrm{perc}}\right)$ be the scaled adjacency matrix of $\mathcal{G}_{\mathrm{perc}}\left(\mathcal{G}_{\mathrm{lat}},\left\{p\right\}_{d=1}^{d=D}\right)$ for factor $\gamma=\sum_{d=1}^{D}p_d\left(M_d-1\right)$ with empirical spectral distribution $F_{W}$, and let $\widehat{A}\left(\mathcal{G}_{\mathrm{perc}}\right)=\Delta^{-1}\left(\mathcal{G}_{\mathrm{perc}}\right)A\left(\mathcal{G}_{\mathrm{perc}}\right)$ be the row-normalized adjacency matrix of $\mathcal{G}_{\mathrm{perc}}\left(\mathcal{G}_{\mathrm{lat}},p\right)$  with empirical spectral distribution $F_{\widehat{A}}$, where $\Delta\left(\mathcal{G}_{\mathrm{perc}}\right)$ is the diagonal matrix of node degress.  Also let $d_{\scriptscriptstyle\mathrm{L}}\left(\cdot,\cdot\right)$ be the L\'{e}vy distance metric.  Assume that all of the lattice dimension sizes increase without bound as $N\rightarrow \infty$.  Then,
\begin{equation}
\lim_{N\rightarrow\infty}{d_{\scriptscriptstyle\mathrm{L}}\left(F_{\sqrt{\gamma}\widehat{A}},F_{\sqrt{\gamma}W}\right)}=0.
\end{equation}
\label{Thm3}
\end{theorem}

\vspace{-20pt}
\section{Conclusion}\label{Conclusion}

This paper analyzed the eigenvalues of the scaled adjacency matrices of a random graph model formed by non-uniform Bernoulli percolation of a $D$-dimensional lattice graph where the link inclusion probability parameter is structured such that it depends only on the lattice dimension index of the link.  Specifically, a deterministic equivalent sequence of distribution functions are computed for the sequence of  empirical spectral distributions of the adjacency matrices using the stochastic canonical equations techniques of Girko in order to capture the asymptotic behavior as the lattice dimension sizes increase without bound.  Theorem \ref{Thm2} derived the form of the solution to an important matrix equation from Theorem \ref{GirkoK01Thm} used to compute the empirical spectral distribution.  Corollary~\ref{Cor1} finds the parameter of this solution form by simultaneously diagonalizing the components of this matrix equation and, thus, computes the Stieltjes transform of the deterministic equivalent distributions.  Simulations demonstrate the results for selected parameters.  Finally, Theorem \ref{Thm3} describes the relationship of the scaled adjacency matrix empirical spectral distribution to the row-normalized adjacency matrix empirical spectral distribution, which asymptotically become very close.  This type of information could be of use in the design of linear shift-invariant filters for graph signal processing, an application of which could, for instance, be accelerated consensus filters for random graphs.  Future efforts will focus on a more precise characterization of asymptotically diminishing density regions, actual use of the eigenvalue information for filter design applications, and extension of this analysis to additional models.

%\newpage
\bibliographystyle{IEEEbib}

\end{document}